\documentclass[11pt]{amsart2000}

\newtheorem{theorem}{Theorem}[section]
\newtheorem{lemma}{Lemma}[section]

\newtheorem{corollary}{Corollary}[theorem]

\newcommand{\nempty}{\not= \emptyset}
\newcommand{\invlim}[2]{\ensuremath{\lim\limits_{\leftarrow}\{#1,#2\}} }

\newcommand{\B}[1]{\ensuremath{\mathbb{#1}}}
\newcommand{\C}[1]{\ensuremath{\mathcal{#1}}}
\newcommand{\chain}[1]{\{\C{#1}_i\}_{i=1}^{\infty}}

\begin{document}
\setlength{\unitlength}{0.01in}
\linethickness{0.01in}
\begin{center}
\begin{picture}(474,66)(0,0) 
\multiput(0,66)(1,0){40}{\line(0,-1){24}}
\multiput(43,65)(1,-1){24}{\line(0,-1){40}}
\multiput(1,39)(1,-1){40}{\line(1,0){24}}
\multiput(70,2)(1,1){24}{\line(0,1){40}}
\multiput(72,0)(1,1){24}{\line(1,0){40}}
\multiput(97,66)(1,0){40}{\line(0,-1){40}} 
\put(143,66){\makebox(0,0)[tl]{\footnotesize Proceedings of the Ninth Prague Topological Symposium}}
\put(143,50){\makebox(0,0)[tl]{\footnotesize Contributed papers from the symposium held in}}
\put(143,34){\makebox(0,0)[tl]{\footnotesize Prague, Czech Republic, August 19--25, 2001}}
\end{picture}
\end{center}
\vspace{0.25in}
\setcounter{page}{253}
\title[Orbits of turning points]{Orbits of turning points for maps of
finite graphs and inverse limit spaces}
\author{Brian Raines}
\address{Mathematical Institute\\
University of Oxford\\
Oxford OX1 3LB\\
United Kingdom}
\email{raines@maths.ox.ac.uk} 
\thanks{Brian Raines,
{\em Orbits of turning points for maps of finite graphs and inverse limit 
spaces},
Proceedings of the Ninth Prague Topological Symposium, (Prague, 2001),
pp.~253--263, Topology Atlas, Toronto, 2002}
\begin{abstract}
In this paper we examine the topology of inverse limit spaces generated by
maps of finite graphs. 
In particular we explore the way in which the structure of the orbits of
the turning points affects the inverse limit. 
We show that if $f$ has finitely many turning points each on a finite
orbit then the inverse limit of $f$ is determined by the number of
elements in the $\omega$-limit set of each turning point. 
We go on to identify the local structure of the inverse limit space at the
points that correspond to points in the $\omega$-limit set of $f$ when the
turning points of $f$ are not necessarily on a finite orbit. 
This leads to a new result regarding inverse limits of maps of the
interval.
\end{abstract}
\maketitle

\section[intro]{Introduction} 

Every one-dimensional continuum is an inverse limit on finite graphs, and
many, though not all, are homeomorphic to an inverse limit on a finite
graph with a single bonding map. These spaces also naturally appear in
dynamical systems. R.F. Williams showed that if a manifold diffeomorphism
$F$ has a one-dimensional hyperbolic attractor $\Lambda$ (with associated
stable manifold structure)then $F$ restricted to $\Lambda$ is
topologically conjugate with the shift homeomorphism on an inverse limit
of a piecewise monotone map $f$ of some finite graph, \cite{williams}, and
Barge and Diamond, \cite{barge&diamond1}, remark that for any map
$f:G\rightarrow G$ of a finite graph there is a homeomorphism
$F:\B{R}^3\rightarrow \B{R}^3$ with an attractor on which $F$ is conjugate
to the shift homeomorphism on $\invlim{G}{f}$. More recently, Anderson
and Putnam, \cite{anderson&putnam}, have shown that the dynamics arising
from a substitution tiling is often conjugate to the action of a shift-map
on an inverse limit of a branched d-manifold. They then demonstrate how
to use knowledge about the inverse limit space to compute the cohomology
and K-theory of a space of tilings. Extending these ideas, Barge,
Jacklitch and Vago, \cite{bargejacklitch&vago}, use inverse limits induced
by certain Markov maps on wedges of circles to analyze one-dimensional
substitution tiling spaces and one-dimensional unstable manifolds of
hyperbolic sets. Many of their results rely on showing that certain pairs
of these inverse limit spaces are not homeomorphic.

It is often quite difficult to distinguish between inverse limit spaces,
even when the dynamics of the bonding maps are very different. Many
papers have been written to this end, \cite{barge&diamond2},
\cite{barge&martin}, \cite{bruin}, \cite{kailhofer}, and \cite{raines}. 
However most of the techniques have been focused on maps of the interval. 
Perhaps one of the easiest ways to decide if two inverse limit spaces are
not homeomorphic is to count their endpoints. Barge and Martin have shown
that the number of endpoints of $\invlim{[0,1]}{f}$ is the same as the
number of elements in the $\omega$-limit set of the turning points of the
bonding map, $f$, when $f$ has a dense orbit and finitely many turning
points, \cite{barge&martin}.

In this paper we distinguish between these inverse limit spaces by showing
that many of the points have neighborhoods homeomorphic to the product of
a zero-dimensional set and $(0,1)$, and we show that the exceptional
points are those that always project onto the $\omega$-limit set of the
turning points. We do this for inverse limits on graphs, but of course,
the result holds for inverse limits on the interval. In the case of the
interval, our theorem is still an extension of Barge and Martin's result,
because there are many bonding maps that give rise to a three-endpoint
indecomposable continua that have more than three points in the
$\omega$-limit set of their turning points. Our theorem can be used to
easily distinguish between these inverse limit spaces.

\section[prelim]{Preliminaries}\label{prelim} 

By a {\em continuum} we mean a compact, connected, metric space, and by a
{\em mapping} we mean a continuous function. We will say a mapping, $f$,
is {\em monotone} on $A$ if, and only if, $f^{-1}(x)$ is connected for all
$x\in A$. The {\em inverse limit} induced by a single bonding map, $f$, on
a continuum $M$ is defined as follows: $$\invlim{M}{f}=\{(x_0,
x_1,\dots)|x_i\in M \quad\mbox{and}\quad f(x_{i+1})=x_i\}.$$ Since $M$ is
metric and $f$ is a mapping, $\invlim{M}{f}$ is a continuum with the
metric: $$d(x,y)=\sum_{i=0}^{\infty}\frac{d_M(x_i-y_i)}{2^i},$$ where
$d_M$ is the metric on $M$ and we assume that $d_M(x,y)<1$ for all $x,y\in
M$. Define the projection maps $\pi_n:\invlim{M}{f}\rightarrow M$ by
$\pi_n(x)=x_n$, where $x=(x_1, x_2, \dots)\in \invlim{M}{f}$.
Also, define the shift homeomorphism 
$h:\invlim{M}{f}\rightarrow \invlim{M}{f}$ by 
$$h(x)=(f(x_0),f(x_1),f(x_2),\dots)=(f(x_0),x_0,x_1,\dots).$$

A {\em linear chaining}, or just {\em chaining}, of a continuum $M$ is a
finite sequence, $L_1,L_2,L_3,\dots ,L_n$ of open subsets of $M$ such that
$L_i$ intersects $L_j$ if and only if $|i-j|<2$. The open sets comprising
the chain are called the {\em links} of the chain. The {\em mesh} of a
chain is the largest of the diameters of its links. A continuum $M$ is
said to be {\em chainable} provided that for each positive number
$\epsilon$ there is a chaining of $M$ with mesh less than $\epsilon$, such
a chain is called an {\em $\epsilon$-chain}. It is a well-known fact that
an inverse limit of chainable continua is a chainable continuum. A {\em
closed chain} is a chain whose links are closed sets and if $i\not= j$,
then $L_i\cap L_j=$Bd$(L_i)\cap $Bd$(L_j)$ if, and only if, $|i-j|<2$. We
lose no generality in assuming that all of the chains in this paper are
taut, i.e. if $L_i\cap L_j=\emptyset$ then $\overline{L_i}\cap
\overline{L_j}=\emptyset$, \cite{nadler}. Notice that if $\C L$ is a taut
chaining of an inverse limit space and $L_i\cap L_j=\emptyset$ then it is
possible to find a positive integer, $q$, large enough so that
$\pi_q(\overline{L_i})\cap \pi_q(\overline{L_j})=\emptyset$, which implies
that $\pi_q(\C L)$ is a chain.

A {\em finite graph}, $G$, is a continuum that can be written as the union
of finitely many arcs any two of which are either disjoint or intersect at
only one of their endpoints. For any finite graph, $G$, there is a finite
set of points called {\em vertices}, $V=\{v_1,v_2,\dots,v_n\}$, and a set
of arcs, $E$, with endpoints from $V$ called {\em edges}, with the
property that if $v_k\in e_{ij}\in E$ then either $k=i$ or $k=j$, and if
two edges meet, they meet only at a single common vertex. For simplicity,
we adopt the convention that, since $e_{ij}=e_{ji}$, if we label an edge
$e_{ij}$ then $i<j$. For every point, $x\in G$, define the {\em degree of
$x$}, deg$(x)$, to be the number of edges in $G$ that have $x$ as an
endpoint. Let $V'\subseteq V$ be the set of all points, $x$, with
deg$(x)\ge 3$.

Let $a,b\in G$. We will denote an arc between $a$ and $b$ by
$\overline{ab}$. Clearly this arc is not uniquely determined. However if
$a,b\in e_{ij}$ then there is a unique arc with endpoints $a$ and $b$ that
is contained in $e_{ij}$. We will denote this arc by $[a,b]$, assuming
that $a$ is closer to $v_i$ in the linear ordering of $e_{ij}$ that has
$v_i$ as its least element, otherwise we denote it $[b,a]$.

We will now extend the idea of linear-chains to graph-chains. Let $n$ be a
positive integer and let $R$ be a relation on $\{1,2\dots n\}$ with the
property that if $(i,j)\in R$ then $i<j$. For every $(i,j)\in R$, let $\C
C_{i,j}=\{ C^{i,j}_1, C^{i,j}_2, \dots, C^{i,j}_{n_{i,j}}\}$ be a taut
chain with the closure of every link of $\C C_{i,j}$ being disjoint from
the closure of every link of $\C C_{k,l}$ whenever $(k,l)\not= (i,j)$,
except $C^{i,j}_1$ which meets every link of the form $C^{i,k}_1$ and
every $C^{m,i}_{n_{m,i}}$ and also except for $C^{i,j}_{n_{i,j}}$ which
meets every link of the form $C^{j,k}_1$ and every $C^{m,j}_{n_{m,j}}$. 
Call $\C C_{i,j}$ an {\em edge-chain}. Let 
$$\C C = \bigcup_{(i,j)\in R} \C C_{i,j}.$$ 
Call $\C C$ a {\em
graph-chain}. Let $G$ be a finite graph with vertex set, $V$, and edge
set $E$. A {\em graph-chaining of $G$} is a graph-chain with
$R=\{(i,j)|e_{ij}\in E\}$ such that each vertex, $v_i$, is only in links
of the form $C^{i,j}_1$ or $C^{m,i}_{n_{m,i}}$, and each edge-chain, $\C
C_{i,j}$, is a chaining of the corresponding edge, $e_{ij}$.

For a given graph-chain $\C C$, call the set $E'=R$ the {\em edge index
set}. For notational convenience we will often denote a graph-chain by
$$\C C=\{C^i_1, C^i_2, \dots C^i_{n_i}|i\in E'\}$$ 
using $i$ to represent an ordered pair in $E'$.

A continuum, $M$, is said to be {\em graph-chainable} provided that for
each positive number $\epsilon$ there exists a graph-chaining of $M$ with
mesh less than $\epsilon$. By a {\em closed graph-chain} we will mean a
graph-chain $\C{C}$ such that every link of $\C{C}$ is closed and if $A$
and $B$ are different links in $\C{C}$ with $A\cap B\nempty$, then $A\cap
B=$Bd$(A)\cap$Bd$(B)$. Notice that this implies that the only point in
common to links of the form $C^{i,j}_1$ and $C^{k,i}_{n_{k,i}}$ is the
vertex $v_i$.

It is easy to show that the inverse limit induced by maps on
graph-chainable continua is itself a graph-chainable continuum.

\section{Markov Graph-Maps}\label{markovmap}

First, we extend the definition of a Markov map of the interval (see
\cite{barge&diamond2} or \cite{holte1})to a Markov map of a graph. Let
$f$ be a mapping of a finite graph, $G$, with vertex set $V$ and edge set
$E$ and define a {\em Markov graph-chaining}, $\C{T}^f$, of $G$ with
respect to $f$ to be a closed graph-chaining of $G$ where, for every $i\in
E'$, $|\C{T}^f_i|=n_i$, $f$ restricted to each link is monotone but not
constant, and for every $i\in E'$ and $k\le n_i$ there exists a subset of
$E'\times \B{N}$, $A_{i,k}$, such that 
$$f(T^i_k)=\bigcup_{(p,r)\in A_{i,k}}T^p_r.$$ 
We will call a set of the form $A_{i,k}$ the {\em index
set} for $(i,k)$ under $f$, and we will call a map that admits such a
Markov graph-chain a {\em Markov graph-map} or simply a {\em Markov map}. 
The endpoints of each link of $\C{T}^f$ determine a Markov partition of
each edge. We define the {\em Markov partition} of the graph to be the
set 
$$B_f = 
\{v_i=c^{i,j}_0<c^{i,j}_1<\cdots <c^{i,j}_{n_{i,j}}=v_j|(i,j)\in E'\}$$ 
where $c^{i,j}_k$ and $c^{i,j}_{k+1}$ are the endpoints of 
$T^{i,j}_{k+1}$. 
Notice that $f(B_f)\subseteq B_f$, $f$ is not constant on 
$[c^{i,j}_{k-1}, c^{i,j}_k]$, and $f$ restricted to each such arc is
monotone.

Also define the set $S_{i,k}\subset E'\times \B{N}$ such that 
$(p,r)\in S_{i,k}$ if and only if 
$[f^{-1}(T^i_k)]^{\circ}\cap T^p_r\nempty$ i.e.\ $S_{i,k}$ is the
collection of indices of links of $\C{T}_f$ that are mapped onto $T^i_k$
by $f$. 
We will call the set $S_{i, k}$ the {\em inverse index set} of $(i,k)$
under $f$.

Let $f:G_1\rightarrow G_2$ be a a map between finite graphs $G_1$ and
$G_2$. Then $x\in G_1$ is a {\em turning point} of $f$ if there is an
arc, $\overline{ab} \subseteq G_1$, containing $x$ in its interior, such
that $f[\overline{ab}]=\overline{zf(x)}$ where $z\in \{f(a), f(b)\}$ and
both of $f|_{A}$ and $f|_{B}$ are monotone, where
$A=\overline{ax}\subseteq \overline{ab}$ and $B=\overline{xb}\subseteq
\overline{ab}$. Denote the set of turning points of $f$ by $P_f$.

Generally there is much freedom in determining the Markov chain; however
we assume that the Markov chains used in this paper are ``natural'' in the
sense that elements of the Markov partition are either vertices, turning
points, or in the orbit of a turning point or vertex. In the next section
we will need to assume that $f^{-1}(x)$ consists of only isolated points,
for all $x\in G$. The next theorem demonstrates that we lose no generality
in assuming this when $f$ is Markov.

\begin{theorem}\label{bigone} 
Let each of $f$ and $g$ be a Markov mapping of $G$, a finite graph, with
associated Markov partitions, 
$B_f=\{c^i_0<c^i_1<\cdots <c^i_{n_i}|i\in E'\}$ and
$B_g=\{d^i_0<d^i_1<\cdots <d^i_{n_i}|i\in E'\}$. 
Suppose that for every $i\in E'$ and $k\le n_i$ there is a $p\in E'$ and
a $r\le n_p$ such that $f(c^i_k)=c^p_r$ if and only if $g(d^i_k)=d^p_r$,
then $\invlim{G}{f}$ is homeomorphic to $\invlim{G}{g}$. 
\end{theorem}

Before presenting the proof of this theorem we will present a few useful
facts about graph-chainable continua.

Let $\C{C}$ be a closed graph-chaining of a continuum, $M$, with edge
index set $R$, and let $\C C'$ be a refinement of $\C C$ with edge index
set $R'$. Let $h$ be a function such that for every link, $C^{'i}_k\in \C
C'$, let $h(i, k)=(p, r)$ if and only if $C^{'i}_k$ is a subset of $C^p_r$
in $\C{C}$. In this case we shall say that $\C{C}'$ {\em follows pattern
$h$} in $\C{C}$. The proof of the next theorem is quite obvious, and so it
has been omitted.

\begin{theorem}\label{smallone}
Let $A$ be a graph-chainable continuum, and let $\chain{C}$ be a sequence
of refining graph-chainings of $A$ such that 
$$\lim\limits_{i\rightarrow \infty} \operatorname{mesh}(\C{C}_i)=0$$ 
and 
$\C{C}_i$ follows
pattern $h_i$ in $\C{C}_{i-1}$. 
If $B$ is also a graph-chainable continuum with a sequence of refining 
graph-chainings, $\chain{D}$, such that 
$$\lim\limits_{i\rightarrow \infty} \operatorname{mesh}(\C{D}_i)=0$$ 
and $\C{D}_i$ follows pattern $h_i$ in $\C{D}_{i-1}$ then $A$ is 
homeomorphic to $B$. 
\end{theorem}

Suppose that $f:G\rightarrow G$ is a Markov mapping of a finite graph,
$G$, with vertex set $V=\{v_1,v_2,\dots,v_n\}$ and edge set $E$, and let
$\C{T}^f$ be a Markov chain of $G$ for $f$ where, for every $i\in E'$,
$|\C{T}^f_i|=n_i$. Suppose that $\C{C}$ is a closed refinement of
$\C{T}^f$ with $|\C{C}_i|=m_i$, such that $\C{C}$ follows pattern $h$ in
$\C{T}^f$. We define the {\em Markov graph-chain function}, $\hat{f}$, on
the elements of $\C{C}$ by the following. (We denote the lexicographical
ordering on $E'$ by $\ll$.)

First let $j$ be the least integer, $k$, such that $(1,k)\in E'$ and define
$$\hat{f}_{p,r}(C^{1,j}_1)=f^{-1}(C^{1,j}_1)\cap T^p_r$$
where $(p,r)\in S_{h((1,j),1)}$. For $(k,\ell)\in E'$ and
$(p,r)\in S_{h((k,\ell), 1)}$ let
$$
\begin{array}{lll}
\hat{f}_{p,r}(C^{k,\ell}_1)&
=&
\displaystyle{
\left[ 
f^{-1}(C^{k,\ell}_1)\cap T^p_r
\right]
-
\left[
\bigcup_{(q,s)\in E', (q,s)\ll (k, \ell)} \hat{f}_{p,r}(C^{q,s}_1)
\right]^{\circ}
.
}
\end{array}
$$

For $(p,r)\in S_{h((1,j), n_{i,j})}$ define
$$
\begin{array}{lll}
\hat{f}_{p,r}(C^{1,j}_{n_{1,j}})&
=&
\displaystyle{
\left[
f^{-1}(C^{1,j}_{n_{1,j}})\cap T^p_r 
\right]
-
\left[
\bigcup_{(j,\ell)\in E'}\hat{f}_{p,r}(C^{j,\ell}_1)
\right]^{\circ}
.
}
\end{array}
$$

For $(k, \ell)\in E'$ and $(p,r)\in S_{h((k, \ell), n_{k, \ell})}$, 
let
$$
\begin{array}{llll}
\multicolumn{4}{l}{
\hat{f}_{p,r}(C^{k, \ell}_{n_{k, \ell}}) =
}
\\
&
&
\multicolumn{2}{l}{
f^{-1}(C^{k, \ell}_{n_{k, \ell}})\cap T^p_r 
}
\\
&
&
&
\displaystyle{
-
\left(
\left[
\bigcup_{(\ell, m)\in E'}\hat{f}_{p,r}(C^{\ell, m}_1)
\right]^{\circ} 
\cup 
\left[ 
\bigcup_{(q,s)\in E', (q,s)\ll (k, \ell)}\hat{f}_{p,r}(C^{q,s}_{n_{q,s}})
\right]^{\circ}
\right)
.
}
\end{array}
$$

Finally for any $(j,k)\in E'$, $1<m<n_{j,k}$ and 
$(p,r)\in S_{h((j,k), m)}$ 
let
$$
\begin{array}{llll}
\hat{f}_{p,r}(C^{j,k}_m)&
=&
\multicolumn{2}{l}{
f^{-1}(C^{j,k}_m)\cap T^p_r
}
\\
&
&
&
- 
\left( 
\left[
\hat{f}_{p,r}(C^{j,k}_{m-1})\right]^{\circ}\cup \left[
\hat{f}_{p,r}(C^{j,k}_1)\right]^{\circ}\cup \left[
\hat{f}_{p,r}(C^{j,k}_{n_{j,k}})\right]^{\circ} 
\right)
.
\end{array}
$$

Define
$$\hat{f}(\C{C})=\left\{ \hat{f}_{p,r}(C^{i,j}_k)|(i,j)\in E',\
k\le n_i\quad\mbox{and}\quad (p,r)\in S_{h((i,j),k)} \right \}.$$

Notice that since $f$ restricted to each link of $\C{T}^f$ is
monotone, $f$ restricted to each link of $\C{C}$ is monotone. So
each element of $\hat{f}(\C{C})$ is connected.

\begin{lemma}\label{lemma1} 
If $\C{C}$ is a closed refinement of $\C{T}^f$ then $\hat{f}(\C{C})$ is a
closed graph-chain and the components of $\hat{f}(\C{C})$ refine 
$\C{T}^f$. 
\end{lemma}

\begin{proof}
For every $i\in E'$, let $m_i=|\C{C}_i|$. 
Every element of $\hat{f}(\C{C})$ is closed and $\hat{f}(\C{C})$ covers 
$G$. 
Suppose now that $C^i_k\cap C^p_r=\emptyset$, but 
$\hat{f}(C^i_k)\cap \hat{f}(C^p_r)\nempty$. 
Let $x\in \hat{f}(C^i_k)\cap \hat{f}(C^p_r)$. 
This implies that $f(x)\in C^i_k\cap C^p_r$, a contradiction. 
So the only elements of $\hat{f}(\C{C})$ which intersect are images of
links of $\C{C}$ which intersected, and since $\C{C}$ is a closed
graph-chaining of $G$, $\hat{f}(\C{C})$ is also a closed graph-chaining of
$G$.
By definition, links of $\hat{f}(\C{C})$ intersect only on their boundary
and the components of $\hat{f}(\C{C})$ refine $\C{T}^f$.
\end{proof}

Now suppose that $g$ is another Markov mapping of $G$ and let
$$\C{S}^g=\{S^i_0, \dots S^i_{n_i}|i\in E'\}$$ 
be the Markov graph-chain associated with $g$ where, for every $i\in E'$, 
$|\C{S}_i|=n_i$. 
Denote the corresponding Markov partition by 
$$B_g=\{d^{i,j}_0<d^{i,j}_1 < \cdots <d^{i,j}_{n_{i,j}}|(i,j)\in E'\}.$$

We are now ready to prove Theorem \ref{bigone}.

\begin{proof}[Proof of \ref{bigone}]
Choose a positive number $\delta_1$ such that if $\C{H}$ is a closed
graph-chaining of $G$ with mesh less than $\delta_1$ which refines
$\C{T}^f$ or $\C{S}^g$ then $\pi_1^{-1}(\C{H})\cap \invlim{G}{f}$ and
$\pi_1^{-1}(\C{H})\cap \invlim{G}{g}$ both have mesh less than
$\frac{1}{2}$. For every $i\in E'$ and $j\le n_i$ , let $Q^i_j$ be a
positive integer such that 
$\delta_1 \cdot Q^i_j > \operatorname{diam}(T^i_j)$ and
$\delta_1 \cdot Q^i_j > \operatorname{diam}(S^i_j)$.

Let $\C{C}_1$ be a closed graph-chaining of $G$ with mesh less than
$\delta_1$ which refines $\C{T}^f$ such that, for every $i\in E'$ and
$j\le n_i$, there are $Q^i_j$ links of $\C{C}_1$ contained in $T^i_j$. 
Let $\C{J}_1$ be defined similarly with respect to $g$ and $\C{S}^g$.

Let $\C{D}_1=\pi_1^{-1}(\C{C}_1)\cap \invlim{G}{f}$ and let
$\C{K}_1=\pi_1^{-1}(\C{J}_1)\cap \invlim{G}{g}$. It is easy to see that
both of $\C{D}_1$ and $\C{K}_1$ are closed graph-chainings of
$\invlim{G}{f}$ and $\invlim{G}{g}$ respectively with mesh less than
$\frac{1}{2}$.

By lemma \ref{lemma1} both of $\hat{f}(\C{C}_1)$ and $\hat{g}(\C{J}_1)$
are closed graph-chainings of $G$ which refine $\C{T}^f$ and $\C{S}^g$
respectively. Let $\delta_2$ be a positive number so that any closed
graph-chaining of $G$, $\C{H}$, with mesh less than $\delta_2$ has
$\pi_2^{-1}(\C{H})\cap \invlim{G}{f}$ and $\pi_2^{-1}(\C{H})\cap
\invlim{G}{g}$ both have mesh less than $\frac{1}{4}$.

By the hypothesis of the theorem, $\hat{f}_{p,r}(C^i_j)$ is defined if and
only if $\hat{g}_{p,r}(J^i_j)$ is defined. Notice that by the construction
of $\C{C}_1$ and $\C{J}_1$, there is a function, $\ell$, such that
$\C{C}_1$ follows pattern $\ell$ in $\C{T}^f$ and $\C{J}_1$ also follows
pattern $\ell$ in $\C{S}^g$. So for every $i\in E'$, $j\le n_i$, and
$(p,r)\in S_{\ell(i,j)}$, let $Q^{p,r}_{i,j}$ be a positive integer such
that 
$Q^{p,r}_{i,j} \cdot \delta_2 > 
\operatorname{diam}(\hat{f}_{p,r}(C^i_j))$ and
$Q^{p,r}_{i,j}\cdot \delta_2 > \operatorname{diam}(\hat{g}_{p,r}(J^i_j))$. 
Let $\C{C}_2$
be a refinement of $\hat{f}(\C{C}_1)$ such that there are $Q^{p,r}_{i,j}$
links of $\C{C}_2$ inside each $\hat{f}_{p,r}(C^i_j)$. Let $\C{J}_2$ be a
refinement of $\hat{g}(\C{J}_2)$ defined similarly. Define $\C{D}_2$ to
be $\pi_2^{-1}(\C{C}_2)\cap \invlim{G}{f}$ and define $\C{K}_2$ to be
$\pi_2^{-1}(\C{J}_2)\cap \invlim{G}{g}$.

Notice that if $A$ is a subset of $\hat{f}_{p,r}(C^i_j)$ then
$$
\pi_2^{-1}(A)\cap \invlim{G}{f}\subseteq \pi_1^{-1}(C^i_j)\cap 
\invlim{G}{f},
$$ 
and similarly if $A$ is a subset of $\hat{g}_{p,r}(J^i_j)$ then 
$$
\pi_2^{-1}(A)\cap \invlim{G}{g}\subseteq \pi_1^{-1}(J^i_j)\cap
\invlim{G}{g}.
$$ 
So, since we have exactly $Q^{p,r}_{i,j}$ links of $\C{C}_2$ and $\C{J}_2$
in $\hat{f}_{p,r}(C^i_j)$ and $\hat{g}_{p,r}(J^i_j)$ respectively,
$\C{D}_2$ follows the same pattern, $h_2$, in $\C{D}_1$ that $\C{K}_2$
follows in $\C{K}_1$.

Clearly, chains of $\invlim{G}{f}$ and $\invlim{G}{g}$, $\C{D}_3$ and
$\C{K}_3$ can be constructed such that 
$\operatorname{mesh}(\C{D}_3) < \frac{1}{8}$,
$\operatorname{mesh}(\C{K}_3) < \frac{1}{8}$, 
and both $\C{D}_3$ and $\C{K}_3$ follow pattern $h_3$ in $\C{D}_2$ and 
$\C{K}_2$ respectively.

So it is easy to see that we can build a sequence of refining chainings,
$\chain{D}$, of $\invlim{G}{f}$ such that $\C{D}_i$ follows pattern $h_i$
in $\C{D}_{i-1}$ and 
$$\lim\limits_{i\rightarrow \infty} \operatorname{mesh}(\C{D}_i)=0,$$
and we can build a sequence of refining chainings,
$\chain{K}$, of $\invlim{G}{g}$ such that $\C{K}_i$ follows pattern $h_i$
in $\C{K}_{i-1}$ and 
$$\lim\limits_{i\rightarrow \infty} \operatorname{mesh}(\C{K}_i)=0.$$
Thus, by Theorem \ref{smallone}, $\invlim{G}{f}$ is homeomorphic to 
$\invlim{G}{g}$.
\end{proof}

This theorem provides some justification for the assumption that we make
in the next section that $f^{-1}(x)$ is completely disconnected for all
$x\in G$. It shows that for a bonding map with finitely many turning
points each on a finite orbit that we lose no generality in making this
assumption. It also has an interesting, and immediate, corollary which is
an extension of a theorem of Holte (\cite{holte1}, Theorem 3.2).

\begin{corollary} 
Let each of $f$ and $g$ be Markov maps of the interval with associated
Markov partitions, $B_f=\{0=c_0<c_1< \cdots <c_n=1\}$, and
$B_g=\{0=d_0<d_1<\cdots <d_n=1\}$. 
Suppose that $f(c_i)=c_j$ if, and only if, $g(d_i)=d_j$, then 
$\invlim{[0,1}{f}$ is homeomorphic to $\invlim{[0,1]}{g}$.
\end{corollary}

This corollary extends Holte's theorem in the sense that the Markov
partitions for $f$ and $g$ do not need to be the same set of points in the
interval. So it allows for not only eliminating flat spots, but also
dramatically changing slopes and shifting turning points around.

\section{Inverse Limit Spaces}

In this section we will assume that $f:G\to G$ is continuous, has finitely
many turning points and, for every $x\in G$, $f^{-1}(x)$ is completely
disconnected. We will also assume that for every $y\in G$ there exists a
positive integer $n$ such that $f^{-n}(y)$ consists of more than one point
and if $C\subseteq G$ is connected then 
$\operatorname{diam}(C') \le \operatorname{diam}(C)$ for every
component, $C'$ of $f^{-1}(C)$.

Let $E_G\subseteq G$ be the set of endpoints of $G$. The { \em
$\omega$-limit set of a point $x$ under a mapping $f$}, $\omega_f(x)$ or
simply $\omega(x)$, is given by $$\omega(x)=\bigcap_{N\in \B N}\overline{
\{f^n(x)|n\ge N\} },$$ and the $\omega$-limit set of a set, $X$, is given
by $\omega(X)=\bigcup_{x\in X} \omega(x).$

Denote the elements of $P_f$, the turning points of $f$, in the edge
$e_{i,j}$ by $t^{i,j}_{1,0}<t^{i,j}_{2,0}\cdots <t^{i,j}_{m,0}$, where $<$
is the linear ordering of the arc $e_{i,j}$ that has $v_i$ as its least
element, and for every $t^{i,j}_{k,0}$ denote the orbit of $t^{i,j}_{k,0}$
by 
$$
\mbox{orb}(t^{i,j}_{k,0}) = 
\{t^{i,j}_{k,\ell}=f^{\ell}(t^{i,j}_{k,0})|\ell \in \B N\}.
$$ 
Let
$$
\mbox{orb}(P_f) = \bigcup_{t^{i,j}_{k,0}\in P_f} \mbox{orb}(t^{i,j}_{k,0}).
$$

\begin{theorem}
Let $x\in \invlim{G}{f}$ have the property that for some $N\in \B N$, if
$n\ge N$ then deg$(x_n)\in \{0,2\}$. 
Then there is a positive number $\epsilon$ and a zero-dimensional set $S$
with the property that $B_{\epsilon}(x)$ is homeomorphic to $(0,1)\times
S$ if, and only if:
\begin{itemize}
\item[(i)]
there is a positive number, $\delta$, and a positive integer, $s$, such
that $B_{\delta}(x_s)\cap f^p(E_G)=\emptyset$ for every $p\ge 0$, and
\item[(ii)]
there exists a positive integer, $m$, such that $x_m\not \in \omega
(P_f)$.
\end{itemize}
\end{theorem}

\begin{proof}
Let $N\in \B N$ such that if $n\ge N$ then deg$(x_n)<3$, and let $m \ge N$
be such that $x_m\not \in \omega(P_f)$. Since $f[\omega(P_f)]\subseteq
\omega(P_f)$, if $n\ge m$ then $x_n \not \in \omega(P_f)$. Let $n \ge m$. 
Since $x_n \not \in \omega(P_f)$, $x_n \not \in \omega(t^{i,j}_{k,0})$ for
every $t^{i,j}_{k,0}\in P_f$. By the definition of the $\omega$-limit
set, for every $t^{i,j}_{k,0}\in P_f$, there is a positive integer,
$p^{i,j}_{k,0}$ with the property that $x_n\not \in \overline{
\{t^{i,j}_{k,r}|r \ge p^{i,j}_{k,0}\} }$. Since there are only finitely
many turning points for $f$, there exists a positive integer, $q$, such
that $x_n\not \in \overline{ \{t^{i,j}_{k,r}|r\ge q\} }$ for all
$t^{i,j}_{k,0} \in P_f$. So $x_{n+q+1}$ is not in $\overline{
\{t^{i,j}_{k,r}|r\in \B N\} }$ for all $t^{i,j}_{k,0}\in P_f$. Let
$p=n+q+1$. There is a positive number, $\lambda_p< \delta$, such that
$$B_{\lambda_p}(x_p)\cap \left[ \bigcup_{t^{i,j}_{k,0}\in P_f} \{
t^{i,j}_{k,r}|r\in \B N \} \right]=\emptyset,$$ and if $y\in
B_{\lambda_p}(x_p)$ then deg$(y)$ is either $0$ or $2$. So
$B_{\lambda_p}(x_p)$ is homeomorphic to $(0,1)$, and if $q\in \B N$ and
$A$ is any component of $f^{-q}[ B_{\lambda_p}(x_p)]$ then $A$ is
homeomorphic to $(0,1)$. Let $\epsilon$ be a positive number such that
$\pi_p [ B_{\epsilon}(x)]\subseteq B_{\lambda_p}(x_p)$. Let $A=\pi_p [
B_{\epsilon}(x)]$. Let $A_1, A_2, \dots A_n$ be the components of
$f^{-1}(A)$. For each $i\le n$, $f|_{A_i}$ is monotone and $A_i\cap
f^m(E)=\emptyset$ for all $m\ge 0$. For each $i\le n$, let $A^i_1, A^i_2,
\dots A^i_{n_i}$ be the components of $f^{-1}(A_i)$. Again, for each
$j\le n_i$, $f|_{A^i_j}$ is monotone and $A^i_j\cap f^m(E)=\emptyset$ for
all $m\ge 0$. Assuming that $A^{i,j,\dots p,t}_k$ is a component of
$f^{-1}[A^{i,j,\dots p}_t]$, define $A^{i,j,\dots k}_1, A^{i,j\dots
k}_2\dots A^{i,j,\dots k}_{n_{i,j,\dots k}}$ to be the components of
$f^{-1}[A^{i,j,\dots t}_k]$.

Let $S$ be a collection of sequences of positive integers such that
$\langle y_i\rangle \in S$ if, and only if, $A_{y_1}$ is defined and for
every other $i\in \B N$, $A^{y_1, \dots y_i}_{y_{i+1}}$ is defined.
Clearly $S$ is zero-dimensional. Each sequence, $\langle y_i\rangle $ in
$S$ defines a sequence of open arcs, $A_{y_1}, A^{y_1}_{y_2}, \dots $,
with the property that $f$ maps $A^{y_1,y_2,\dots y_i}_{y_{i+1}}$ onto
$A^{y_1,y_2\dots y_{i-1}}_{y_i}$ monotonically. So
$\invlim{A^{y_1,y_2,\dots y_i}_{y_{i+1}}}{f}$ is homeomorphic to $(0,1)$,
and 
$$
B_{\epsilon}(x) = 
\bigcup_{\langle y_i\rangle \in S}
\invlim{A^{y_1,y_2,\dots y_i}_{y_{i+1}}}{f}
$$ 
is homeomorphic to $S\times (0,1)$.

Now assume that $x$ does not satisfy either (i) or (ii) in the theorem. 
First let $\epsilon$ be a positive number, and consider the
$\epsilon$-neighborhood around $x$, $B_{\epsilon}(x)$. Let $n$ be a
positive integer and let $\gamma$ be a positive number such that
$\pi^{-1}_n[B_{\gamma}(x_n)]\subseteq B_{\epsilon}(x)$, $B_{\gamma}(x_n)$
meets $P_f$ at at most one point, and if $y\in B_{\gamma}(x_n)$ then
deg$(y)<3$.

If for every positive number $\delta<\gamma$ and every positive integer
$M$, there exists a positive integer, $m>M$, such that every image under
$f^{-m}$ of $B_{\delta}(x_n)$ meets the set of endpoints for $G$, $E_G$,
then clearly for every positive number, $\lambda$, we can $\lambda$-chain
$\invlim{G}{f}$ with a linear subchain that starts at $x$. Thus there is
a chainable endcontinuum in $G$ having $x$ as an endpoint, and
$B_{\epsilon}(x)$ cannot be homeomorphic to $(0,1)\times S$ where $S$ is a
zero-dimensional set.

Instead, now suppose that $\gamma$ is small enough and $n$ is large enough
so that $f^{-m}(B_{\gamma}(x_n))$ misses $E_G$ for every positive integer
$m$. Let $A=B_{\gamma}(x_n)$, and as above, enumerate the components of
the preimages of $A$, $A_1, A_2, \dots A_n$. Continuing as previously,
let $S$ be the set of sequences of positive integers, $\langle
y_i\rangle$, where $\langle y_i\rangle \in S$ if, and only if, for every
positive integer $i$, $A^{y_1,\dots y_i}_{y_{i+1}}$ is a component of
$f^{-1}(A^{y_1, \dots y_{i-1}}_{y_i})$. Let $t^{i,j}_{k,0}\in P_f$ with
$x_n\in \omega(t^{i,j}_{k,0})$. Then for infinitely many positive
integers, $m$, $t^{i,j}_{k,m}\in A$. So there are infinitely many
connected subsets of $G$, $A^{y_1\dots y_r}_{y_{r+1}}$ that contain
$t^{i,j}_{k,0}$. If there exists one of these subsets, $A^{y_1,\dots
y_r}_{y_{r+1}}$ that only meet $B_f$ at the singleton $t^{i,j}_{k,0}$
then, since these components do not contain any endpoints or vertices of
$G$, every component of the preimages of $A^{y_1, \dots y_r}_{y_{r+1}}$ is
mapped with a single fold across $A^{y_1, \dots y_r}_{y_{r+1}}$. Thus
$B_{\epsilon}(x)$ contains a subspace homeomorphic to a neighborhood of
$(0,1)$ in the $\sin(1/x)$-continuum, and it cannot be homeomorphic to
$(0,1)\times T$, where $T$ is a zero-dimensional set. If instead, for
some $A^{y_1,\dots, y_r}_{y_{r+1}}$ we have $A^{y_1,\dots
y_r}_{y_{r+1}}\cap B_f$ is a finite set then clearly we can restrict the
size of $A^{y_1, \dots, y_r}_{y_{r+1}}$ in order to make the subset meet
$B_f$ at a singleton and produce a similar subspace.

So suppose that each subset, $A^{y_1,\dots y_r}_{y_{r+1}}$ that contains a
turning point meets $B_f$ on an infinite set. Pick one of these subsets
and call it $A_1$. Let $A_2$ be a connected subset of $G$ such that
$f^{n_1}(A_2)=A_1$ and $A_2\cap P_f\nempty$. Given $A_i$ define $A_{i+1}$
to be a connected subset of $G$ with the property that
$f^{n_i}(A_{i+1})=A_i$ and $A_{i+1}$ meets $P_f$. Then since $A$ is small
enough to not meet $P_f$ at two points and since all of these subsets are
preimages of $A$, they must all meet $P_f$ at a single point. Similarly,
they must all miss $V'$ and the set of endpoints of $G$. Thus, for every
positive integer, $i$, $f|_{A_i}$ is a two-pass map, and
$f^{n_i}|_{A_{i+1}}$ is at least a two-pass map. Thus
$\invlim{A_i}{f|_{A_i}}$ is an indecomposable subcontinuum, and
$B_{\epsilon}(x)$ contains an indecomposable subcontinuum. Hence,
$B_{\epsilon}(x)$ is not homeomorphic to $(0,1)\times T$ where $T$ is a
zero-dimensional set.
\end{proof}

\begin{corollary}
Suppose that $f$ and $g$ are maps of $[0,1]$ with the properties listed
above. 
Further suppose that $|\omega(P_f)|=n$ and $|\omega(P_g)|=m$. 
If $n\not= m$ then $\invlim{[0,1]}{f}$ is not homeomorphic to
$\invlim{[0,1]}{g}$.
\end{corollary}

\begin{proof} 
Notice that for every point in $\omega(P_f)$ there is a point in the
inverse limit that either is an endpoint or is in a neighborhood
homeomorphic to a neighborhood of $(0,1)$ in the 
$\sin(1/x)$-continuum. 
Also notice that any other point in the inverse limit has a neighborhood
homeomorphic to the product of $(0,1)$ with a zero-dimensional set, $S$. 
These properties are preserved by homeomorphism, and so any space
homeomorphic to it must have the same number of points in the 
$\omega$-limit set of its turning points.
\end{proof}

\providecommand{\bysame}{\leavevmode\hbox to3em{\hrulefill}\thinspace}
\providecommand{\MR}{\relax\ifhmode\unskip\space\fi MR }
\providecommand{\MRhref}[2]{%
  \href{http://www.ams.org/mathscinet-getitem?mr=#1}{#2}
}
\providecommand{\href}[2]{#2}

\end{document}